\documentclass[a4paper]{jpconf}

\usepackage{graphicx}

\usepackage{eufrak}

\usepackage{amsbsy,amssymb,amscd,amsfonts,latexsym,amstext,delarray,amsmath,color,caption}
\usepackage{hyperref}
\usepackage{tikz}
\usetikzlibrary{matrix,arrows}
\usepackage{graphicx}

\usepackage{hyphenat}

\usepackage{soul}

\usepackage{graphicx,relsize}

\usepackage{cite}

\usepackage{mathrsfs}
\usepackage[errorshow]{tracefnt}
\usepackage{fancybox}
\usepackage{amsfonts}

\usepackage{mathdots}
\usepackage{multirow}
\usepackage{array}
\usepackage{mathtools}
\usepackage{soul}
\usepackage{pict2e}
\usepackage{mathabx}
\DeclareMathAlphabet{\mathpzc}{OT1}{pzc}{m}{it}

\usepackage{float}

\usepackage{xypic}
\usepackage{color}
\usepackage[
color,matrix,arrow]{xy}
\usepackage{setspace}

\usepackage{enumitem}
\setitemize{label=\usebeamerfont*{itemize item}%
  \usebeamercolor[fg]{itemize item}
  \usebeamertemplate{itemize item}}

\usepackage{graphicx}

\let\latexcirc=\circ
\newcommand{\ccirc}{\mathbin{\mathchoice
  {\xcirc\scriptstyle}
  {\xcirc\scriptstyle}
  {\xcirc\scriptscriptstyle}
  {\xcirc\scriptscriptstyle}
}}
\newcommand{\xcirc}[1]{\vcenter{\hbox{$#1\latexcirc$}}}
\let\circ\ccirc

\newcommand{\scr}{\mathscr}

\def\fg{\mathfrak{g}}

\def\small#1{{\hbox{$#1$}}}

\def\boxit#1{\vbox{\hrule\hbox{\vrule\kern3pt
             \vbox{\kern3pt#1\kern3pt}\kern3pt\vrule}\hrule}}

\newcommand{\beq}{\begin{equation}}
\newcommand{\beqn}{\begin{equation*}}
\newcommand{\eeq}{\end{equation}}
\newcommand{\eeqn}{\end{equation*}}
\newcommand{\beqa}{\begin{eqnarray}}
\newcommand{\beqan}{\begin{eqnarray*}}
\newcommand{\eeqa}{\end{eqnarray}}
\newcommand{\eeqan}{\end{eqnarray*}}
\newcommand{\bdm}{\begin{displaymath}}
\newcommand{\edm}{\end{displaymath}}

\newcommand{\ba}{\begin{array}}
\newcommand{\ea}{\end{array}}

\renewcommand{\theenumii}{\alph{enumii}}

\newcommand\nn{\nonumber}

\newcommand\benu{\begin{enumerate}}
\newcommand\eenu{\end{enumerate}}
\newcommand\bit{\begin{itemize}}
\newcommand\eit{\end{itemize}}

\def\der'{\mathfrak{der}'\,}
\def\der{\mathfrak{der}\,}
\def\str'{\mathfrak{str}'\,}
\def\str{\mathfrak{str}\,}

\def\so{\mathfrak{so}}

\def\sl{\mathfrak{sl}}

\newcommand{\de}{\delta}

\renewcommand\arraystretch{1.5}

\newcommand{\adiagram}{
\begin{picture}(250,30)(45,-10)
\put(48,-10){${\scriptstyle{0}}$}
\put(88,-10){${\scriptstyle{1}}$}
\put(128,-10){${\scriptstyle{2}}$}
\put(242,-10){${\scriptstyle{n-2}}$}
\put(282,-10){${\scriptstyle{n-1}}$}
\thicklines
\put(50,10){\line(1,1){3.5}}
\put(50,10){\line(-1,1){3.5}}
\put(50,10){\line(1,-1){3.5}}
\put(50,10){\line(-1,-1){3.5}}
\multiput(50,10)(40,0){3}{\circle{10}}
\multiput(250,10)(40,0){2}{\circle{10}}
\multiput(55,10)(40,0){2}{\line(1,0){30}}
\put(135,10){\line(1,0){20}}
\put(165,10){\line(1,0){10}}
\put(185,10){\line(1,0){10}}
\put(205,10){\line(1,0){10}}
\put(225,10){\line(1,0){20}}
\put(255,10){\line(1,0){30}}
\end{picture}
}

\RequirePackage{hyperref}

\begin{document}
\title{Generators and relations for (generalised) Cartan type superalgebras}

\author{Lisa Carbone${}^1$, Martin Cederwall${}^2$ and Jakob Palmkvist${}^2$}

\address{${}^1$Department of Mathematics, Rutgers University,
  110 Frelinghuysen Rd, Piscataway, NJ~05584, USA}
\address{${}^2$Department of Physics, Chalmers University of Technology, SE-412 96 G\"oteborg, Sweden}

\ead{lisa.carbone@rutgers.edu, martin.cederwall@chalmers.se, jakob.palmkvist@chalmers.se}

\begin{abstract}
In Kac's classification of finite-dimensional Lie superalgebras, the contragredient ones can be constructed from Dynkin diagrams similar to those of the simple finite-dimensional Lie algebras, but with additional types of nodes. For example, $A(n-1,0) = \mathfrak{sl}(1|n)$ can be constructed by adding a ``gray'' node to the Dynkin diagram of $A_{n-1} = \mathfrak{sl}(n)$, corresponding to an odd null root. The Cartan superalgebras constitute a different class, where the simplest example is $W(n)$, the derivation algebra of the Grassmann algebra on $n$ generators.
Here we present a novel construction of $W(n)$, from the same Dynkin diagram as $A(n-1,0)$, but with additional generators and relations.
\end{abstract}

This talk, given by JP at ``The 32nd International Colloquium on Group Theoretical Methods in Physics (Group32)'' in Prague, 9--13 July, 2018, is based on \cite{Carbone:2018xqq}, where more details and references can be found.

In attempts to understand the origin of the duality symmetries appearing in supergravity theories, the Lie algebras 
$\fg$
describing the symmetries have been extended to infinite-dimensional Lie superalgebras.
These extensions include Borcherds superalgebras \cite{HenryLabordere:2002dk}, here denoted by $\scr B(\fg)$, as well
as tensor hierarchy algebras \cite{Palmkvist:2013vya}, here denoted by $W(\fg)$ and $S(\fg)$.

The construction of tensor hierarchy algebras in \cite{Palmkvist:2013vya} was only applicable for finite-dimensional $\fg$
and thus in particular not to the cases
$\fg=E_r$ for $r\geq 9$. Here we solve this problem by a new construction with generators and relations. We focus on the
case $\fg=A_{n-1}$ 
since $W(\fg)$ and $S(\fg)$ 
then turn out to be finite-dimensional and well known as {\it Cartan type} superalgebras. In this sense, the general tensor hierarchy algebras
$W(\fg)$ and $S(\fg)$ are {\it generalised} Cartan type superalgebras.

We consider algebras over an algebraically closed field $\mathbb{K}$ of characteristic zero.
A {\it superalgebra} $G$ is an algebra with a $\mathbb{Z}_2$-grading, which means that it can be decomposed into a direct sum 
$G=G_{(0)}\oplus G_{(1)}$ of an even subalgebra $G_{(0)}$ and an odd subspace $G_{(1)}$,
such that $G_{(i)}G_{(j)} = G_{(i+j)}$ where $i,j \in \mathbb{Z}_2$. (This means that $G_{(1)}$ does not close under the product and is thus
not a subalgebra.)
In a {\it Lie} superalgebra, the product is a bracket that satisfies the identities
\begin{align}
[x,y]&=-(-1)^{|x||y|}[y,x]\,, \label{supersymmetry}\\
[x,[y,z]]&=[[x,y],z]+(-1)^{|x||y|}[y,[x,z]]\,, \label{jacobi-identity}
\end{align}
where $|x|=0$ if $x\in G_{(0)}$ and $|x|=1$ if $x\in G_{(1)}$.

A $\mathbb{Z}$-{\it grading} of the Lie superalgebra $G$ is a
decomposition of $G$  into a direct sum of subspaces $G_i$ for all
integers $i$,  called {\it levels}, such that $[G_i,G_j]\subseteq
G_{i+j}$. The $\mathbb{Z}$-grading is said to be {\it consistent} if $G_i \subseteq G_{(i)}$ for all levels $i \in \mathbb{Z}$,
that is, if odd elements appear at odd levels and even elements at even levels.

One way to obtain a Lie superalgebra if we only have a superalgebra $G$ to start with (not necessarily a Lie superalgebra) is to consider all even and odd derivations of it. Let $|d|$ be equal to either 1 or 0. A linear map $d: G \to G$ satisfying 
\begin{align}
d(xy)&=d(x)y+(-1)^{|d||x|}xd(y)
\end{align}
is an odd derivation if $|d|=1$ and an even derivation if $|d|=0$. All even and odd derivations of $G$ span a vector space
which, together with the commutator as the bracket,
\begin{align}
[c,d]=c \circ d - (-1)^{|c||d|} d \circ c,
\end{align}
forms a Lie superalgebra. This is the {\it derivation algebra} of $G$, denoted ${\rm der}\, G$.

The Grassmann algebra $\Lambda(n)$ is a basic example of a superalgebra that is not a Lie superalgebra. It is
the associative superalgebra generated by $n$ odd elements
$\theta^0,\theta^1,\ldots,\theta^{n-1}$,
modulo the relations $\theta^a\theta^b = -
\theta^b\theta^a$, where $a,b=0,1,\ldots,n-1$. It is spanned by monomials $\theta^{a_1} \cdots \theta^{a_p}$, where $0 \leq p \leq n$,
which are fully antisymmetric in the upper indices,
\begin{align}
\theta^{a_1} \cdots \theta^{a_p}= \theta^{[a_1} \cdots \theta^{a_p]}\,.
\end{align}
The derivation algebra of $\Lambda(n)$ has a basis consisting of elements
\begin{align}
K^{a_1 \cdots a_p}{}_b = \theta^{a_1} \cdots \theta^{a_p}
\frac{\partial}{\partial\theta^b}
\end{align}
acting on a monomial $\theta^{c_1}\cdots\theta^{c_q}$ by a contraction,
\begin{align} \label{contraction}
K^{a_1 \cdots a_p}{}_b \quad : \quad \theta^{c_1}\cdots\theta^{c_q}
\mapsto q\,\delta_b{}^{[c_1}\theta^{|a_1} \cdots \theta^{a_p|}\theta^{c_2}\cdots \theta^{c_q]}\,.
\end{align}
This Lie superalgebra ${\rm der}\,\Lambda(n)$ is also denoted by $W(n)$.
It is easy to see that it has a consistent $\mathbb{Z}$-grading
where the subspace at level $-p+1$
has a
basis of elements $K^{a_1 \cdots a_p}{}_b$ which are fully antisymmetric in
the upper indices, and thus there are no elements at level $-p+1$ for $p>n$ (or $p<0$).
Negative and positive levels
are here reversed compared to the usual conventions.
\begin{align} \label{W-basis}
\begin{array}{c|c}
\text{level}&\text{basis}\\
\hline
 1 & K_a\\
 0 & K^{a}{}_b \\
 -1 & K^{ab}{}_c\\
 \vdots &\vdots\\
 -\,n+1 & K^{a_1\cdots a_n}{}_b
\end{array}
\end{align}

In the classification of simple finite-dimensional Lie superalgebras, $W(n)$ appear as Lie superalgebras of {\it Cartan type}. Such Lie superalgebras are distinguished from the {\it classical} Lie superalgebras, which are further divided into {\it basic} and {\it strange} ones \cite{Kac77B}.

The basic Lie superalgebras are finite-dimensional cases of {\it contragredient} Lie superalgebras \cite{Kac77B}, which means that they can be constructed
from a (generalised) Cartan matrix, or from a Dynkin diagram encoding the same information as the matrix.
For a general contragredient Lie superalgebra, the only condition on the Cartan matrix is that it be a square matrix $B_{ab}$
with values in $\mathbb{K}$, where the index set labelling rows and columns
is $\mathbb{Z}_2$-graded, that is, the disjoint union of an odd and an even subset. Here we restrict to the cases where 
$B_{ab}$ ($a,b=0,1,2\ldots,r$) is obtained by adding a row and a column to the Cartan matrix $A_{ij}$
($i,j=1,2\ldots,r$) of a Kac--Moody algebra $\fg$ of rank $r$, 
such that
\begin{align}
B_{ij}&=A_{ij}\,,&B_{0a}&=B_{a0}=\begin{cases}-1 &(a=1)\,, \\ \phantom{-}0  & 
(a\neq1)\,, \end{cases}
\end{align}
and such that $\{0\}$ and $\{1,2,\ldots,r\}$ are the odd and even subsets, respectively, of the 
index set $\{0,1,2,\ldots,r\}$. Furthermore, for simplicity we assume that 
$A_{ij}$ is symmetric (implying that $B_{ab}$ is symmetric as well) and that both $A_{ij}$ and $B_{ab}$ are non-degenerate.
To this Cartan matrix we associate a $\mathbb{Z}_2$-graded set $M=M_{(0)} \cup M_{(1)}$ of generators, 
\begin{align} \label{generators}
M_{(0)}&= \{e_i,f_i,h_a\}\,& M_{(1)}&=\{e_0,f_{0}\}\,,
\end{align}
where $i=1,2,\ldots,r$ and $a=0,1,2,\ldots,r$.
Let $\scr B(\fg)$ be the Lie superalgebra generated by the set $M$ modulo the Chevalley--Serre relations
\begin{align*} 
[h_a,e_b]&=B_{ab}e_b\,, & [h_a,f_b]&=-B_{ab}f_b\,,& [e_a,f_b]&=\delta_{ab}h_b\,,
\end{align*}
\begin{align} \label{serre-rel0}
({\rm ad}\ e_a)^{1-B_{ab}}(e_b)&=({\rm ad}\ f_a)^{1-B_{ab}}(f_b)=0\,.
\end{align}
Then $\scr B(\fg)$ is the contragredient Lie superalgebra constructed from the Cartan matrix $B_{ab}$.
With the restrictions on $B_{ab}$ here, $\scr B(\fg)$ is not only a contragredient Lie superalgebra
but also a {\it Borcherds} superalgebra \cite{Borcherds,Ray}.
    
\begin{figure}
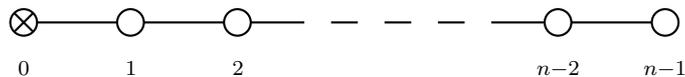

\begin{center}
\adiagram 
\end{center}
\caption{\label{superADynkin}The Dynkin diagram of $\scr B(A_{n-1})=A(n-1,0)$.}
\end{figure}

Let us now apply the above construction to the
case of the (finite-dimensional) Kac--Moody algebra $\fg=A_{n-1}$ (thus $r=n-1$) 
with a Cartan matrix $A_{ij}$ where row and column 1 correspond to one of the end nodes in the Dynkin diagram.    
Then we get the Cartan matrix $B_{ab}$ of $\scr B(\fg)$ given in (\ref{cartan-A-case}).      
We can associate a Dynkin diagram to it, given in Figure \ref{superADynkin}, where row and column 0 correspond to the ``gray'' node.
\begin{align} \label{cartan-A-case}
B_{ab}=\begin{pmatrix}
0 & -1 & 0 &\cdots& 0 & 0\\
-1 & 2 & -1 &\cdots& 0 & 0\\
0 & -1 & 2 &\cdots& 0 & 0\\
\vdots &\vdots &\vdots & \ddots & \vdots & \vdots\\ 
0 & 0& 0 & \cdots & 2 & -1\\
0 & 0& 0 & \cdots & -1 & 2
\end{pmatrix}\,
\end{align}
The resulting Lie superalgebra $\scr B(A_{n-1})$ is $A(n-1,0)$, one of the basic Lie superalgebras in the classification \cite{Kac77B}.

If we put $e_0$ and $f_0$ at level 1 and level $-1$, respectively, and all the other generators at level 0, then we get a consistent $\mathbb{Z}$-grading
of $\scr B(\fg)$. In the case $\fg=A_{n-1}$ there are no more levels,
this $\mathbb{Z}$-grading of $\scr B(A_{n-1})$ is a
3-grading:
\begin{align} \label{3gradingofA}
A(n-1,0) = \mathscr B_{-1} \oplus \mathscr B_{0} \oplus \mathscr B_{1}\,.
\end{align}
The subalgebra $\scr B_0$ is $\mathfrak{sl}(n)\oplus\mathbb{K}=\mathfrak{gl}(n)$,
and the basis elements can thus be written as $\mathfrak{gl}(n)$
tensors:
\begin{align} \label{B-basis}
\begin{array}{c|c|c}
\text{generators}&\text{level}&\text{basis}\\
\hline
e_0 & 1 & E_a\\
e_i,f_i,h_a & 0 & G^{a}{}_b \\
f_0 & -1 & F^a
\end{array}
\end{align}

The commutation relations are
\begin{align}
[G^a{}_b,G^c{}_d]&=\de_b{}^cG^a{}_d-\de_d{}^aG^c{}_b\,,  & [E_a,F^b]&=-G^b{}_a +\de_a{}^b G\,,\nn
\end{align}
\begin{align}
[G^a{}_b,F^c]&= \de_b{}^c F^a\,,& [G^a{}_b,E_c]&= -\de_c{}^a E_b\,, & [E_a,E_b]&=[F^a,F^b]=0\,,
\label{A(n-1,0)commrel}
\end{align}
where $G=\sum_{a=0}^{n-1}G^a{}_a$.
With the identifications
\begin{align*}
  e_0 &= E_0\,, & f_0 &= F^0\,,
  & h_0&=G^1{}_1 + G^2{}_2 + \cdots + G^{n-1}{}_{n-1} = G-G^0{}_0\,,
\end{align*} 
\begin{align} \label{chev2}
e_i &= G^{i-1}{}_{i}\,, & f_i &= G^{i}{}_{i-1}\,, & h_i&=G^{i-1}{}_{i-1}-G^{i}{}_{i}\,,
\end{align}
the commutation relations (\ref{A(n-1,0)commrel}) follow from the
Chevalley--Serre relations
(\ref{serre-rel0}). 
    
Let us now compare the basis (\ref{B-basis}) of $A(n-1,0)$ to the basis (\ref{W-basis}) of $W(n)$.
Level 1 and 0 have the same index structure in $W(n)$ as in $A(n-1,0)$. If we consider level $-1$ in $W(n)$,
the tensors can be decomposed into a traceless part
and the trace, obtained by contracting the lower index with one of the upper indices. If we take the subalgebra of $W(n)$ generated by level 1 and only the traceless part of level $-1$, then we get another Lie superalgebra of Cartan type, denoted $S(n)$, with traceless tensors all the way down to level $-n+2$. If we instead take the subalgebra generated by level 1 and only the trace part of level $-1$, then there will be no lower levels, and we get $A(n-1,0)$.
Thus $A(n-1,0)$ is a subalgebra of $W(n)$, which can be constructed from the generators (\ref{generators}) and the relations (\ref{serre-rel0}).
The question arises whether we can obtain not only this subalgebra, but the whole of $W(n)$ by extending the set of generators and relations.

We extend the set $M=M_{(0)}
\cup M{}_{(1)}$ of generators to $M'=M_{(0)}
\cup M'{}_{(1)}$, where $f_0$ is replaced by $r$ generators $f_{0a}$.
\begin{align} \label{extendedgenerators}
M_{(0)}&= \{e_i,f_i,h_a\}\,,& M'{}_{(1)}&=\{e_0,f_{0a}\,|\,a \neq 1\}\,.
\end{align}
Henceforth,
whenever $f_{0a}$ appears we assume $a=0,2,3,\ldots,r$ (with $r=n-1$ in the case of $\fg=A_{n-1}$),
and whenever $f_{a}$ appears we assume $a \neq 0$. As we will see, the new generator $f_{00}$ corresponds to the old $f_0$.
Identifying them with each other, $M'$ is indeed an extension of $M$.

We let $\widetilde W(\fg)$ be the Lie superalgebra generated by the set $M'$ modulo the relations
(\ref{serre-rel0}) and
the additional relations
\begin{align*}
[e_0,f_{0a}]&=h_a\,,& [h_a,f_{0b}]&=-B_{a0}f_{0b}\,,&[e_1,f_{0a}]&=0\,,
\end{align*}
\begin{align*}
[e_a,[e_a,e_{0b}]]=[f_a,[f_a,f_{0b}]]=0\,,
\end{align*}
\begin{align} \label{extendedrelations}
i,j=2,3,\ldots,\,r \ \ \Rightarrow\ \ [e_i,[f_j,f_{0a}]]=\delta_{ij}B_{aj}f_{0j}\,,
\end{align}
where $a=0,1,2,\ldots,r$.

In the same way as for $\scr B(\fg)$ we get a consistent $\mathbb{Z}$-grading of $\widetilde W(\fg)$
if we put $e_0$ at level 1, all $f_{0a}$ at level $-1$, and all the other generators at level 0.
The sum of the subspaces at level $\pm1$ and 0 constitute the {\it local part} of $\widetilde W(\fg)$.
We then define $W(\fg)$ as $W(\fg)=\widetilde W(\fg)/J$, where $J$ is the 
maximal ideal of $\widetilde W(\fg)$ intersecting the local part
trivially. 

The following theorem summarises the main results of \cite{Carbone:2018xqq} (where the proof can be found).\\\\
\noindent
{\bf Theorem.} {\em The Lie superalgebra ${W}(A_{n-1})=\widetilde{W}(A_{n-1})/J$ is 
  isomorphic to $W(n)$. The ideal $J$ of $\widetilde{W}(A_{n-1})$ is generated by the relations
\begin{align} \label{J-relationer}
[f_{0a},f_{0b}]&=[f_{0i},[f_{0j},f_1]]=[(f_{02}-f_{00}),[f_{0j},f_1]]=0\,,
\end{align}
where $i,j=3,\ldots,n-1$. Thus
$W(n)$ has
    generators (\ref{extendedgenerators}) 
        and relations
(\ref{serre-rel0}), (\ref{extendedrelations}) and (\ref{J-relationer}).}\\\\
\noindent
By removing $h_0$ and $f_{00}$ from the set $M'$ of generators we get a subalgebra $S(\fg)$ of $W(\fg)$. In the case $\fg=A_{n-1}$ we have $S(\fg)=S(n)$.
In general, for finite-dimensional $\fg$ this definition of $S(\fg)$ agrees with the definition of the corresponding tensor hierarchy algebra in \cite{Palmkvist:2013vya}. In cases other than $\fg=A_{n-1}$ we do not know whether the relations (\ref{J-relationer}) generate the whole ideal $J$ or if additional relations are needed.

We conclude this talk with an overview of the cases where $\fg$ belongs to the $A$, $D$ or $E$ series of Kac--Moody algebras, with node 1 being the node to which another node is connected when going to the next algebra in the series. Another Lie superalgebra of Cartan type, $H(2r)$, appears as $S(D_r)$.
\renewcommand{\sl}{\mathfrak{sl}}
\renewcommand{\so}{\mathfrak{so}}
\newcommand{\pytteliteta}{\scalebox{0.7} 
{
\begin{picture}(100,20)(-2.5,-2.5)
\thicklines
\multiput(0,0)(15,0){2}{\circle{5}}
\multiput(45,0)(15,0){4}{\circle{5}}
\put(0,0){\line(1,1){1.8}}
\put(0,0){\line(-1,1){1.8}}
\put(0,0){\line(1,-1){1.8}}
\put(0,0){\line(-1,-1){1.8}}
\put(2.5,0){\line(1,0){10}}
\put(17.5,0){\line(1,0){5}}
\put(25.5,0){\line(1,0){3}}
\put(31.5,0){\line(1,0){3}}
\put(37.5,0){\line(1,0){5}}
\multiput(47.5,0)(15,0){3}{\line(1,0){10}}
\put(15,0){\circle{5}}
\end{picture}
}}
\newcommand{\pyttelitetd}{\scalebox{0.7}{
\begin{picture}(100,20)(-2.5,-2.5)
\thicklines
\multiput(0,0)(15,0){2}{\circle{5}}
\multiput(45,0)(15,0){4}{\circle{5}}
\put(0,0){\line(1,1){1.8}}
\put(0,0){\line(-1,1){1.8}}
\put(0,0){\line(1,-1){1.8}}
\put(0,0){\line(-1,-1){1.8}}
\put(2.5,0){\line(1,0){10}}
\put(17.5,0){\line(1,0){5}}
\put(25.5,0){\line(1,0){3}}
\put(31.5,0){\line(1,0){3}}
\put(37.5,0){\line(1,0){5}}
\multiput(47.5,0)(15,0){3}{\line(1,0){10}}
\put(75,2.5){\line(0,1){10}}
\put(75,15){\circle{5}}
\put(15,0){\circle{5}}
\end{picture}
}}
\newcommand{\pyttelitete}{\scalebox{0.7}{
\begin{picture}(90,20)(-2.5,-2.5)
\thicklines
\multiput(0,0)(15,0){2}{\circle{5}}
\multiput(45,0)(15,0){4}{\circle{5}}
\put(0,0){\line(1,1){1.8}}
\put(0,0){\line(-1,1){1.8}}
\put(0,0){\line(1,-1){1.8}}
\put(0,0){\line(-1,-1){1.8}}
\put(2.5,0){\line(1,0){10}}
\put(17.5,0){\line(1,0){5}}
\put(25.5,0){\line(1,0){3}}
\put(31.5,0){\line(1,0){3}}
\put(37.5,0){\line(1,0){5}}
\multiput(47.5,0)(15,0){3}{\line(1,0){10}}
\put(60,2.5){\line(0,1){10}}
\put(60,15){\circle{5}}
\put(15,0){\circle{5}}
\end{picture}
}}
\setlength{\arraycolsep}{6pt}
\renewcommand{\arraystretch}{2}
\begin{align*}
\begin{array}{c|c|c|c}
& \pytteliteta & \pyttelitetd & \pyttelitete\\
\hline
\mathfrak{g} & A_{n-1}=\sl(n) & D_{r}=\so(2r) &  E_r\\
\hline
\mathscr B(\fg) & A(n-1,0)
=\sl(n\,|\,1) & D(d,1)=\mathfrak{osp}(2r\,|\,2) & \text{infinite-dimensional}\\
\hline
S(\fg) & S(n) & H(2r) & \text{infinite-dimensional}\\
\hline
W(\fg) & W(n) & \text{infinite-dimensional} & \text{infinite-dimensional}
\end{array}
\end{align*}

In applications to extended geometry, the $A$, $D$ and $E$ cases above correspond to ordinary, double and exceptional geometry, respectively \cite{Cederwall:2017fjm}.
In cases where so-called ancillary transformations are absent, the Borcherds superalgebras
can be used to derive expressions
for the generalised diffeomorphisms \cite{Palmkvist:2015dea} and furthermore an $L_\infty$ algebra encoding their gauge structure \cite{Cederwall:2018aab,martin-prag}.
When ancillary transformations are present it seems that the Borcherds superalgebra needs to be replaced by a tensor hierarchy algebra \cite{LinftyTHA1,LinftyTHA2}.

\end{document}